\documentclass[12pt]{article}
 
\usepackage[margin=1in]{geometry} 
\usepackage{amsmath,amsthm,amssymb}

\usepackage[utf8x]{inputenc}    
\usepackage[T1]{fontenc}
\usepackage[english]{babel}

\usepackage{color}

\newtheorem*{conj*}{Conjetura}

\theoremstyle{definition}

\newtheorem*{definicao*}{Defini\c c\~ao}

\begin{document}

\title{\bf{The fractional Kontorovich-Lebedev transform}}

\author{Semyon  YAKUBOVICH}
\maketitle

\markboth{\rm \centerline{ Semyon   YAKUBOVICH}}{}
\markright{\rm \centerline{Fractional Kontorovich-Lebedev Transform}}

\begin{abstract} {\noindent We introduce the so-called fractional analog of the Kontorovich-Lebedev transform, modifying its kernel, which is the Macdonald function $K_\nu(z)$.  Properties of this kernel as well as mapping properties of the new index transform are investigated. The inversion formula in suitable spaces of functions is established.   }

\end{abstract}
\vspace{4mm}



\vspace{4mm}

\section{Introduction and preliminary results}

As it is known [3],  the modified Bessel function or Macdonald function $K_\nu(z)$ can be defined by the integral representation

$$K_{\nu}(z)= \int_{0}^\infty e^{-z\cosh u} \cosh(\nu u) du,\quad {\rm Re} (z) >0,\  \nu \in \mathbb{C}.\eqno(1.1)$$
The pure imaginary $\nu=i\tau,\ \tau \in \mathbb{R}$ and positive $x$ correspond to the kernel $K_{i\tau}(x)$ of the reciprocal Kontorovich-Lebedev transforms 

$$(Ff) (\tau)=\int_0^\infty K_{i\tau}(x) f(x) dx,\eqno(1.2)$$

$$f(x)= {2\over \pi^2 x} \int_0^\infty \tau\sinh(\pi\tau) K_{i\tau}(x) (F f)(\tau) d\tau\eqno(1.3)$$
which generate the class of the index transforms [4], where the integration is realized with respect to the argument and the index (a parameter) of the kernel function. The function $K_\nu(x)$ is a fundamental solution of the Bessel homogeneous second order differential equation

$$x^2 y^{\prime\prime} (x) + x y^{\prime} (x) - (x^2+\nu^2) y(x)= 0.\eqno(1.4)$$
It means that $K_\nu(x)$ is the eigenfunction of the differential operator

 $${\cal A}_x \equiv x^2- x{d\over dx} x {d\over dx},\eqno(1.5)$$
i.e. this yields
$${\cal A}_x \ K_{\nu}(x)= - \nu^2 K_{\nu}(x).\eqno(1.6)$$
Moreover, $K_{\nu}(x)$  remains bounded as $x$ tends to
infinity on the real positive line. It has the asymptotic behaviour [3]
$$ K_\nu(x) = \left( \frac{\pi}{2x} \right)^{1/2} e^{-x} [1+
O(1/x)], \qquad x \to +\infty,\eqno(1.7)$$
and near the origin
$$ K_\nu(x) = O\left ( x^{-|{\rm Re}\nu|}\right), \ x \to 0,\eqno(1.8)$$
$$K_0(x) = -\log x + O(1), \ x \to 0. \eqno(1.9)$$

Our goal in this paper is to modify (1.1), involving a positive parameter $\alpha$ and considering the following function

$$K_{\alpha,\nu}(z)= \int_{0}^\infty e^{-z\cosh^\alpha u} \cosh(\nu u) du,\quad {\rm Re} (z) >0,\  \alpha > 0, \nu \in \mathbb{C}.\eqno(1.10)$$
Particular cases comprise 

$$K_{1,\nu}(z) = K_\nu(z),\quad\quad  K_{2,\nu}(z) = {1\over 2}\  e^{-z/2} K_{\nu/2}\left({z\over 2}\right).$$
Accordingly (see (1.2)), we will investigate the following fractional analog of the Kontorovich-Lebedev transform

$$(F_\alpha f)(\tau)=\int_0^\infty K_{\alpha,i\tau}(x) f(x)  dx,\eqno(1.11)$$
establishing in the sequel its mapping properties and proving the corresponding inversion formula.  Taking into account a fundamental role of the Kontorovich-Lebedev transform in the theory of the index transforms [4], this investigation, hopefully,  will give a departure point to establish fractional analogs of the other index transformations.

\subsection{Properties of the function $ K_{\alpha,\nu}(z)$}

We begin, employing in (1.10) the Mellin-Barnes representation for the exponential function (see Entry 8.4.3.1 in [1, Vol. III]), to write

$$ K_{\alpha,\nu}(z)= {1\over 2\pi i} \int_{0}^\infty  \cosh(\nu u) \int_{\gamma-i\infty}^{\gamma+i\infty} \Gamma (s) \left[z\cosh^\alpha u\right]^{-s} ds du,\eqno(1.12)$$
where $\Gamma(s)$ is the Euler gamma-function [1, Vol. III], $\gamma >0,\  0 \le |\arg(z)| < \pi/2$. Hence, if $\left|{\rm Re} (\nu)\right| < \alpha \gamma$, we find the estimate

$$ \int_{0}^\infty  \left| \cosh(\nu u)\right| \int_{\gamma-i\infty}^{\gamma+i\infty} \left| \Gamma (s) \left[z\cosh^\alpha u\right]^{-s} ds\right|  du$$

$$\le |z|^{-\gamma}   \int_{\gamma-i\infty}^{\gamma+i\infty} \left| \Gamma (s) \right| e^{{\rm Im} (s) \arg (z) } |ds|  \int_{0}^\infty  {\cosh({\rm Re}(\nu) u) \over \cosh^{\alpha\gamma} u} \ du  < \infty.$$ 
Therefore, the interchange of the order of integration in (1.12) is permitted, and, calculating the inner integral via Entry 2.4.4.4. in [1, Vol. I], we arrive at the following Mellin-Barnes representation 

$$ K_{\alpha,\nu}(z)= {1\over 8\pi i}  \int_{\gamma-i\infty}^{\gamma+i\infty}   {\Gamma (s)\over \Gamma(\alpha s)} \Gamma\left({\alpha s+\nu\over 2}\right) \Gamma\left({\alpha s-\nu\over 2}\right)  \left({z\over 2^\alpha}\right)^{-s} ds.\eqno(1.13) $$
The latter integral can be expressed in terms of the Fox $H$-function (cf. Entry 8.4.51.11 in [1, Vol. III]), namely,

$$ K_{\alpha,\nu}(z)= {1\over 4} \mathop{H_{1,3}^{3,0}}\left({z\over 2^\alpha} ; \  {(0,\alpha) \atop (0,1), (\nu/2, \alpha/2), (-\nu/2, \alpha/2)}\right).$$

{\bf Theorem 1}. {\it The function $K_{\alpha,\nu}(z)$ satisfies the following integro-differential equation

$$\left( \alpha^2 z{d\over dz} z {d\over dz} \right) \left[ u(z)- I_ -^{2/\alpha} u (z)\right] = \nu^2  u(z)  - \alpha z {d\over dz}  I_ -^{2/\alpha} u (z),\eqno(1.14)$$
where

$$I_ -^{\beta} u (z) = z^\beta \int_1^\infty (x-1)^{\beta-1} u(zx) dx$$
is the right-handed Riemann-Liouville fractional integral $[3]$.}

\begin{proof} Indeed,  the Stirling asymptotic formula for the gamma function  permits to differentiate  consecutively in (1.13) with respect to $z,\ |z| \ge r_0 >0$ under the integral sign. Hence via properties for the gamma function and simple substitutions we derive

$$\left( \alpha^2 z{d\over dz} z {d\over dz} - \nu^2\right)  K_{\alpha,\nu}(z) =   {1\over 2\pi i}  \int_{\gamma-i\infty}^{\gamma+i\infty}  2^{\alpha s}  {\Gamma (s)\over \Gamma(\alpha s)} \Gamma\left({\alpha s+\nu\over 2}+1\right) \Gamma\left({\alpha s-\nu\over 2}+1\right)  z^{-s} ds $$

$$=   {z^{2/\alpha} \over 8\pi i}  \int_{\gamma+2/\alpha-i\infty}^{\gamma+ 2/\alpha+i\infty}  2^{\alpha s} (\alpha s-2)(\alpha s-1)   {\Gamma (s- 2/\alpha)\over \Gamma(\alpha s)} \Gamma\left({\alpha s+\nu\over 2}\right) \Gamma\left({\alpha s-\nu\over 2}\right)  z^{-s} ds $$

$$=   \alpha^2  z{d\over dz} z^{1+1/\alpha}  {d\over dz}    {z^{1/\alpha} \over 8\pi i}   \int_{\gamma+2/\alpha-i\infty}^{\gamma+ 2/\alpha+i\infty}  2^{\alpha s}  {\Gamma (s- 2/\alpha)\over \Gamma(\alpha s)} \Gamma\left({\alpha s+\nu\over 2}\right) \Gamma\left({\alpha s-\nu\over 2}\right)  z^{-s} ds. $$
Taking into account a simple beta-integral 

$$ {\Gamma (s- 2/\alpha)\over \Gamma( s)} = {1\over \Gamma(2/\alpha) }\int_1^\infty (x-1)^{2/\alpha-1} x^{-s} dx, $$
where ${\rm Re}(s)= \gamma+ 2/\alpha,\ \gamma >0$, the previous equalities  and (1.13) imply 

$$\left( \alpha^2 z{d\over dz} z {d\over dz} - \nu^2\right)  K_{\alpha,\nu}(z) =  {\alpha^2 z \over \Gamma(2/\alpha) }  {d\over dz} z^{1+1/\alpha}  {d\over dz}  z^{1/\alpha}  \int_1^\infty (x-1)^{2/\alpha-1} K_{\alpha,\nu}(z x) dx$$

$$=     {\alpha z (\alpha-1)\over \Gamma(2/\alpha) }  {d\over dz}  z^{2/\alpha}   \int_1^\infty (x-1)^{2/\alpha-1} K_{\alpha,\nu}(z x) dx+ {\alpha^2 z^2 \over \Gamma(2/\alpha) }    {d^2\over dz^2}   z^{2/\alpha}  \int_1^\infty (x-1)^{2/\alpha-1} K_{\alpha,\nu}(z x) dx$$

$$=  \left( \alpha^2 z{d\over dz} z {d\over dz} - \alpha z {d\over dz}  \right) I_ -^{2/\alpha} K_{\alpha,\nu} (z)$$
which means (1.14).

\end{proof}

{\bf Remark 1}. For real positive $z=x$ and $\alpha=1$ we use the index law for fractional integrals [3] to get immediately (1.6). 

{\bf Corollary 1}.  {\it The following equality takes place}

$$\left( \alpha^2 z{d\over dz} z {d\over dz} \right)  \left[ K_{\alpha,\nu}(z) - {1\over 2\pi} \int_{-\infty}^\infty {\tau  K_{\alpha,\nu+i\tau}(z) \over \sinh(\pi\tau/2)}  d\tau\right] $$

$$ = \nu^2  K_{\alpha,\nu}(z) -   {\alpha z\over 2\pi} {d\over dz} \int_{-\infty}^\infty {\tau  K_{\alpha,\nu+i\tau}(z) \over \sinh(\pi\tau/2)}  d\tau.\eqno(1.15)$$

\begin{proof} In fact, recalling (1.10) and the Fubini theorem to interchange the order of integration, we have from the proof of Theorem 1

$$\left( \alpha^2 z{d\over dz} z {d\over dz} - \nu^2\right)  K_{\alpha,\nu}(z) =  {\alpha^2 z \over \Gamma(2/\alpha) }  {d\over dz} z^{1+1/\alpha}  {d\over dz}  z^{1/\alpha}  \int_1^\infty (x-1)^{2/\alpha-1} K_{\alpha,\nu}(z x) dx$$

$$=  \alpha^2 z {d\over dz} z^{1+1/\alpha}  {d\over dz}  z^{- 1/\alpha}  \int_0^\infty   {e^{-z\cosh^\alpha u}\over \cosh^2 u}  \cosh(\nu u)  du $$

$$=  \alpha z \left( \alpha {d\over dz} z {d\over dz} - {d\over dz}\right) \int_0^\infty   {e^{-z\cosh^\alpha u}\over \cosh^2 u}  \cosh(\nu u)  du, $$
i.e.

$$\left( \alpha^2 z{d\over dz} z {d\over dz} \right)  \left[ K_{\alpha,\nu}(z) - \int_0^\infty   {e^{-z\cosh^\alpha u}\over \cosh^2 u}  \cosh(\nu u)  du\right]  = \nu^2  K_{\alpha,\nu}(z) $$

$$-   \alpha z {d\over dz} \int_0^\infty   {e^{-z\cosh^\alpha u}\over \cosh^2 u}  \cosh(\nu u)  du.$$
Meanwhile, the Parseval equality for the Fourier transform, (1.10) and Entry 2.5.46.6 in [1, Vol. I] yield

$$\int_0^\infty   {e^{-z\cosh^\alpha u}\over \cosh^2 u}  \cosh(\nu u)  du = {1\over 2\pi} \int_{-\infty}^\infty {\tau  K_{\alpha,\nu+i\tau}(z) \over \sinh(\pi\tau/2)}  d\tau$$
which gives (1.15).

\end{proof}

{\bf Remark 2}. Letting $\alpha = 1$ in (1.15), we find via (1.6) the following identity for the Macdonald function

$$  K_{\nu}(z) =   {1\over 2\pi} {d^2\over dz^2} \int_{-\infty}^\infty {\tau  K_{\nu+i\tau}(z) \over \sinh(\pi\tau/2)}  d\tau.$$
Further, for two positive parameters $\alpha, \beta$ we follow similar arguments to write from (1.12)

$$ K_{\alpha+\beta,\nu}(z)= {1\over 2\pi i} \int_{0}^\infty  \cosh(\nu u) \int_{\gamma-i\infty}^{\gamma+i\infty} \Gamma (s) \left[z\cosh^{\alpha+\beta} u\right]^{-s} ds du$$

$$= {1\over 2\pi^2 i} \int_{0}^\infty  \int_{\gamma-i\infty}^{\gamma+i\infty} {\Gamma (s) 2^{\beta s-1} \over \Gamma(\beta s)}  \Gamma\left({\beta s+ i\tau\over 2}\right) \Gamma\left({\beta s- i\tau\over 2}\right)  z^{-s} \int_{0}^\infty  {\cosh(\nu u)\cos(\tau u) \over   \cosh^{\alpha s} u } du ds d\tau$$

$$= {1\over 2\pi^2 i}  \int_{\gamma-i\infty}^{\gamma+i\infty} {\Gamma (s) 2^{(\alpha+\beta) s-3} \over \Gamma(\alpha s) \Gamma(\beta s)}  z^{-s} \int_{-\infty}^\infty  \Gamma\left({\beta s\over 2} +i\tau\right) \Gamma\left({\beta s\over 2} -i\tau\right) $$

$$\times  \Gamma\left({\alpha s+ \nu\over 2} +i\tau\right) \Gamma\left({\alpha s- \nu\over 2}-i\tau\right)  d\tau ds.\eqno(1.16)$$
The integral with respect to $\tau$ can be calculated explicitly in terms of residues at simple left-hand poles of gamma functions and the Gauss hypergeometric function [1, Vol. III]. Precisely, we have

$${1\over 2\pi} \int_{-\infty}^\infty  \Gamma\left({\beta s\over 2} +i\tau\right) \Gamma\left({\beta s\over 2} -i\tau\right)  \Gamma\left({\alpha s+ \nu\over 2} +i\tau\right) \Gamma\left({\alpha s- \nu\over 2}-i\tau\right)  d\tau$$

$$= \sum_{n=0}^\infty {(-1)^n\over n!} \bigg[ \Gamma(\beta s+n)  \Gamma\left({(\alpha +\beta) s- \nu\over 2} + n\right) \Gamma\left({(\alpha -\beta) s+ \nu\over 2} - n\right)\bigg.$$

$$\bigg. +  \Gamma(\alpha s+n)  \Gamma\left({(\alpha +\beta) s+\nu\over 2} + n\right) \Gamma\left({(\beta-\alpha) s- \nu\over 2} - n\right)\bigg]$$

$$=  \Gamma\left(\beta s\right) \Gamma\left({(\alpha +\beta) s- \nu\over  2}\right)\Gamma\left({(\alpha -\beta) s+ \nu\over  2}\right) \ {}_2F_1 \left( \beta s,\  {(\alpha +\beta) s- \nu\over 2};\ 1- {(\alpha -\beta) s+ \nu\over 2}; 1 \right)$$

$$+  \Gamma\left(\alpha s\right) \Gamma\left({(\alpha +\beta) s+\nu\over  2}\right) \Gamma\left({(\beta-\alpha) s- \nu\over 2} \right) \ {}_2F_1 \left( \alpha s,\  {(\alpha +\beta) s+ \nu\over 2};\ 1- {(\beta-\alpha) s- \nu\over 2}; 1 \right).$$
Finally, under the condition $\gamma(\alpha+\beta) < 1$ we employ the value of the Gauss hypergeometric function at the unity to obtain after simplification

$${1\over 2\pi} \int_{-\infty}^\infty  \Gamma\left({\beta s\over 2} +i\tau\right) \Gamma\left({\beta s\over 2} -i\tau\right)  \Gamma\left({\alpha s+ \nu\over 2} +i\tau\right) \Gamma\left({\alpha s- \nu\over 2}-i\tau\right)  d\tau$$

$$=   \frac{\pi \Gamma(1- (\alpha+\beta) s) }{\sin(\pi ((\alpha -\beta) s+ \nu)/2) }\bigg[ \frac{ \Gamma\left(\beta s\right)  \Gamma\left(((\alpha +\beta) s- \nu)/2\right)}{\Gamma\left(1- ((\alpha +\beta) s+ \nu)/2\right)  \Gamma(1- \alpha s)}\bigg.$$

$$\bigg. -  \frac{\Gamma\left(\alpha s\right) \Gamma\left(((\alpha +\beta) s+ \nu)/2\right)}{\Gamma\left(1- ((\alpha +\beta) s- \nu)/2\right)  \Gamma(1- \beta s)}\bigg]$$

$$=   \frac{\pi  \Gamma\left(((\alpha +\beta) s- \nu)/2\right) \Gamma\left(((\alpha +\beta) s+ \nu)/2\right)}{\Gamma((\alpha+\beta) s)   \sin(\pi(\alpha+\beta)s)\sin(\pi ((\alpha -\beta) s+ \nu)/2) }$$

$$\times \bigg[ \frac{ \Gamma\left(\beta s\right) \sin(\pi ((\alpha +\beta) s+ \nu)/2)}{ \Gamma(1- \alpha s)} -  \frac{\Gamma\left(\alpha s\right) \sin(\pi ((\alpha +\beta) s- \nu)/2) }{  \Gamma(1- \beta s)}\bigg].$$
Substituting  this expression in (1.16), we derive the following Mellin-Barnes representation for $ K_{\alpha+\beta,\nu}(z)$

$$ K_{\alpha+\beta,\nu}(z) = {1\over 8\pi i}  \int_{\gamma-i\infty}^{\gamma+i\infty}  2^{(\alpha+\beta) s}\   { \Gamma (s) \over \Gamma((\alpha+\beta) s) }  \Gamma\left({(\alpha+\beta) s+\nu\over 2}\right) \Gamma\left({(\alpha+\beta) s-\nu\over 2}\right) z^{-s}$$

$$\times  \frac{ \sin(\pi\alpha s) \sin(\pi ((\alpha +\beta) s+ \nu)/2) -  \sin(\pi\beta s) \sin(\pi ((\alpha +\beta) s- \nu)/2) }{\sin(\pi(\alpha+\beta)s)\sin(\pi ((\alpha -\beta) s+ \nu)/2) } ds$$
under conditions $\gamma >0,\  0 \le |\arg(z)| < \pi/2,\  \alpha, \beta > 0,\   \left|{\rm Re} (\nu)\right| < (\alpha+\beta) \gamma < 1.$ But
$$ \frac{ \sin(\pi\alpha s) \sin(\pi ((\alpha +\beta) s+ \nu)/2) -  \sin(\pi\beta s) \sin(\pi ((\alpha +\beta) s- \nu)/2) }{\sin(\pi(\alpha+\beta)s)\sin(\pi ((\alpha -\beta) s+ \nu)/2) }  = 1.$$
So,  this gives (1.13) for $K_{\alpha+\beta,\nu}(z) $, omitting the condition $(\alpha+\beta) \gamma < 1$.

On the other hand, returning to (1.16) and employing the duplication formula for gamma function, we have

$$ K_{\alpha+\beta,\nu}(z)= {1\over 8\pi\sqrt\pi  i}  \int_{-\infty}^\infty  \int_{\gamma-i\infty}^{\gamma+i\infty} {\Gamma (s)\   \Gamma\left(\beta s/ 2 +i\tau\right)\ \Gamma\left(\beta s/ 2 -i\tau\right)  \over \Gamma(\alpha s)\ \Gamma(\beta s/2) \ \Gamma((\beta s+1)/2)}  $$

$$\times  \Gamma\left({\alpha s+ \nu\over 2} +i\tau\right) \Gamma\left({\alpha s- \nu\over 2}-i\tau\right) \left({z\over 2^\alpha}\right)^{-s}  ds d\tau.\eqno(1.17) $$
Meanwhile, Entry 8.4.41.4 in [1, Vol. III] gives the following Mellin transform of the associated Legendre function of the first kind

$$\int_0^1 (1-x)^{-1/4} P^{1/2}_{-1/2+i\tau}\left({2\over x}-1\right) x^{(\beta s-3)/2} dx =   {\Gamma\left(\beta s/ 2 +i\tau\right)\ \Gamma\left(\beta s/ 2 -i\tau\right)  \over  \Gamma(\beta s/2) \ \Gamma((\beta s+1)/2)},\ {\rm Re}(s) > 0. $$
Making a simple substitution, the latter integral reads

$$\int_0^\infty {x^{-1/4}\over (x+1)^{(1+ 2\beta s)/4}}\  P^{1/2}_{-1/2+i\tau}\left(2x+1\right) dx =   {\Gamma\left(\beta s/ 2 +i\tau\right)\ \Gamma\left(\beta s/ 2 -i\tau\right)  \over  \Gamma(\beta s/2) \ \Gamma((\beta s+1)/2)}.\eqno(1.18) $$
Using the uniform bound (3.95) in [4] for the associated Legendre function

$$\left| P^{1/2}_{-1/2+i\tau}\left(2x+1\right) \right| \le C e^{(\pi - 2\delta)|\tau|} x^{1/4} (x+1)^{-1/4},$$
where $x >0,\ \tau \in \mathbb{R},\ \delta \in [0, \pi/2),\ C >0$ is an absolute constant, it yields the estimate
$$\int_0^\infty \left|{x^{-1/4}\over (x+1)^{(1+ 2\beta s)/4}}\  P^{1/2}_{-1/2+i\tau}\left(2x+1\right)\right|dx \le C  e^{(\pi - 2\delta)|\tau|} \int_0^\infty {dx\over (x+1)^{(1+\beta \gamma)/2}},$$
where the latter integral converges under the condition $\gamma > 1/\beta$. Hence (1.17) takes the form of the iterated integral

$$ K_{\alpha+\beta,\nu}(z)= {1\over 8\pi\sqrt\pi  i}  \int_{-\infty}^\infty  \int_{\gamma-i\infty}^{\gamma+i\infty} {\Gamma (s) \over \Gamma(\alpha s)} \Gamma\left({\alpha s+ \nu\over 2} +i\tau\right) \Gamma\left({\alpha s- \nu\over 2}-i\tau\right) \left({z\over 2^\alpha}\right)^{-s}  $$

$$\times   \int_0^\infty {x^{-1/4}\over (x+1)^{(1+ 2\beta s)/4}}\  P^{1/2}_{-1/2+i\tau}\left(2x+1\right)   dx ds d\tau.\eqno(1.19)$$
Our goal is to interchange the order of integration in (1.19) and appeal then to (1.13).  Indeed,  letting $ \left|{\rm Re} (\nu)\right| < \alpha \gamma,\  \alpha, \gamma >1/\beta,\ \beta >0$, we use an analog of the inequality (1.100) in [4] for the modified Bessel function and Entry 8.4.23.1 in [1, Vol. III] to find 

$$\left| 2^{\alpha s-2} \Gamma\left({\alpha s+ \nu\over 2} +i\tau\right) \Gamma\left({\alpha s- \nu\over 2}-i\tau\right) \right| = \left| \int_0^\infty K_{\nu+2i\tau}(x) x^{\alpha s-1} dx \right|$$

$$\le \int_0^\infty \left| K_{\nu+2i\tau}(x) \right| x^{\alpha \gamma-1} dx \le  e^{-\delta \left| {\rm Im}(\nu)+ 2\tau\right|} \int_0^\infty  K_{{\rm Re}(\nu)}(x) x^{\alpha \gamma-1} dx,$$ 
where  $\delta \in [0, \pi/2)$.  Thus, choosing $0 < \alpha <1,\ \delta \in (\pi/4, \pi/2)$ and taking  $ 0 \le |\arg(z)| < (1-\alpha) \pi/2$, we deduce the estimate

$$\int_{-\infty}^\infty  \int_{\gamma-i\infty}^{\gamma+i\infty} \left|{\Gamma (s) \over \Gamma(\alpha s)} \Gamma\left({\alpha s+ \nu\over 2} +i\tau\right) \Gamma\left({\alpha s- \nu\over 2}-i\tau\right) \left({z\over 2^\alpha}\right)^{-s} \right| $$

$$\times   \int_0^\infty \left|{x^{-1/4}\over (x+1)^{(1+ 2\beta s)/4}}\  P^{1/2}_{-1/2+i\tau}\left(2x+1\right)   dx ds\right| d\tau $$

$$\le C  |z|^{-\gamma} \int_{-\infty}^\infty  \exp \bigg[\pi  |\tau| - \delta \bigg[ 2 |\tau|+  \left| 2\tau + {\rm Im}(\nu)\right|\bigg]\  \bigg] d\tau \int_{\gamma-i\infty}^{\gamma+i\infty} \left|{\Gamma (s) \over \Gamma(\alpha s)} \right|  e^{{\rm Im} (s) \arg (z) } |ds|$$

$$\times  \int_0^\infty {dx\over (x+1)^{(1+\beta \gamma)/2}}   \int_0^\infty  K_{{\rm Re}(\nu)}(y) y^{\alpha \gamma-1} dy < \infty,$$
where  $\gamma > \max (1/\beta,   \left|{\rm Re} (\nu)\right|/\alpha)$. Consequently,  the Fubini theorem allows to interchange the order of integration in (1.18) and by virtue of (1.13) to complete the proof of the following theorem.

{\bf Theorem 2.} {\it Let $0< \alpha <1, \beta >0, \ \nu,\ z  \in \mathbb{C}, \  0 \le |\arg(z)| < (1-\alpha) \pi/2$.  Then function $K_{\alpha+\beta,\nu}(z)$ has the following  double integral representation}

$$  K_{\alpha+\beta,\nu}(z)= {1\over \sqrt\pi}  \int_{-\infty}^\infty  \int_0^\infty   {P^{1/2}_{-1/2+i\tau}\left(2x+1\right) \over (x(x+1) )^{1/4}}\   K_{\alpha, \nu+2i\tau} \left( z(x+1)^{\beta/2}\right) dx  d\tau.\eqno(1.20)$$

{\bf Remark 3}.  Taking into account the parameter symmetry,   a companion  formula  takes place, accordingly,

  $$  K_{\alpha+\beta,\nu}(z)= {1\over \sqrt\pi}  \int_{-\infty}^\infty  \int_0^\infty   {P^{1/2}_{-1/2+i\tau}\left(2x+1\right) \over (x(x+1) )^{1/4}}\   K_{\beta, \nu+2i\tau} \left( z(x+1)^{\alpha/2}\right) dx  d\tau\eqno(1.21)$$
under conditions $0< \beta <1, \alpha >0, \ \nu,\ z  \in \mathbb{C}, \  0 \le |\arg(z)| < (1-\beta) \pi/2$. 

{\bf Corollary 2}. {\it The modified Bessel function has the following decomposition

 $$  K_{\nu}(z)= {1\over \sqrt\pi}  \int_{-\infty}^\infty  \int_0^\infty   {P^{1/2}_{-1/2+i\tau}\left(2x+1\right) \over (x(x+1) )^{1/4}}\   K_{\alpha, \nu+2i\tau} \left( z(x+1)^{(1-\alpha)/2}\right) dx  d\tau,\eqno(1.22)$$
where $0< \alpha < 1, \nu,\ z  \in \mathbb{C},\   0 \le |\arg(z)| < (1-\alpha) \pi/2$}. 

{\bf Corollary 3}. {\it The  function  $K_{\alpha,\nu}(z)$ admits the integral  representation

 $$  K_{\alpha, \nu}(z)= {\sqrt\pi\over 2}   \int_0^\infty  P^{1/2}_{-1/2+ \nu}\left(2x+1\right) \ {e^{-z (x+1)^{\alpha/2}}\over (x(x+1) )^{1/4}}\  dx,\eqno(1.23)$$
where $\alpha > 0, \ \nu,  z  \in \mathbb{C}: \  {\rm Re} (\nu) < 1/2,\  0 \le |\arg(z)| <  \pi/2$}. 

\begin{proof} Rewriting (1.18) as follows

$$\int_0^\infty {x^{-1/4}\over (x+1)^{(1+ 2\alpha s)/4}}\  P^{1/2}_{-1/2+\nu}\left(2x+1\right) dx =   {\Gamma\left(\alpha s/ 2 + \nu\right)\ \Gamma\left(\alpha s/ 2 - \nu\right)  \over  \Gamma(\alpha s/2) \ \Gamma((\alpha s+1)/2)},$$
where $ {\rm Re}(s) >  2 |{\rm Re} (\nu)|/\alpha$,  we see that equality (1.13) takes the form of the iterated integral

$$ K_{\alpha,\nu}(z)= {1\over 4\sqrt \pi i}  \int_{\gamma-i\infty}^{\gamma+i\infty}   \Gamma (s) z^{-s} \int_0^\infty {x^{-1/4}\over (x+1)^{(1+ 2\alpha s)/4}}\  P^{1/2}_{-1/2+\nu}\left(2x+1\right) dx ds.$$
The interchange of the order of integration is justified by the following estimate (cf. (3.94) in [4])

$$ \int_{\gamma-i\infty}^{\gamma+i\infty}   \left|\Gamma (s) z^{-s} \right| \int_0^\infty \left| {x^{-1/4}\over (x+1)^{(1+ 2\alpha s)/4}}\  P^{1/2}_{-1/2+\nu}\left(2x+1\right) dx ds\right|$$

$$\le  |z|^{-\gamma} \int_{\gamma-i\infty}^{\gamma+i\infty}   \left|\Gamma (s) \right|  e^{{\rm Im} (s) \arg (z) }  |ds|  \int_0^\infty {x^{-1/4}\over (x+1)^{(1+ 2\alpha \gamma)/4}}\  \left| P^{1/2}_{-1/2+\nu}\left(2x+1\right) \right| dx$$

$$\le  |z|^{-\gamma} \left|\cos(\pi\nu)\right| \int_{\gamma-i\infty}^{\gamma+i\infty}   \left|\Gamma (s) \right|  e^{{\rm Im} (s) \arg (z) }  |ds|  \int_0^\infty {dx\over (x+1)^{(1+ \alpha \gamma)/2}} $$

$$\times    \int_0^\infty K_{2 {\rm Re}(\nu)} (y) dy < \infty$$
when $z  \in \mathbb{C},\   0 \le |\arg(z)| <  \pi/2,\ \alpha \gamma > 1, \  |{\rm Re} (\nu)| < 1/2$. Therefore, interchanging the order of integration and calculating the inner Mellin-Barnes integral, we end up with (1.23), completing the proof of Corollary 3.

\end{proof}

\subsection{A multiplicative operator  semigroup}

For positive $z=x >0$ we recall (1.13) and reciprocal formulas for the Mellin transforms [2] to conclude the equality

$$\int_0^\infty   K_{\alpha,\nu}(x) x^{s-1} dx =  {2^{\alpha s-2}\  \Gamma (s)\over \Gamma(\alpha s)} \Gamma\left({\alpha s+\nu\over 2}\right) \Gamma\left({\alpha s-\nu\over 2}\right),\eqno(1.24)$$  
where $\left|{\rm Re} (\nu)\right| < \alpha \gamma$.  Let $0 < \alpha <1$.  Then one can introduce the following auxiliary kernel $\Phi_\alpha(x),\ x >0$, being represented by the absolutely convergent Mellin-Barnes integral

$$\Phi_\alpha(x)= {1\over 2\pi i} \int_{\gamma-i\infty}^{\gamma+i\infty} { \Gamma (s)\over \Gamma(\alpha s)} x^{-s} ds.\eqno(1.25)$$
This function can be easily expressed in terms of the Wright function as a sum of residues in the left-handed  simple poles of the gamma function, namely,

$$\Phi_\alpha(x)=  {1\over \pi}  \sum_{n=0}^\infty  (-1)^{n+1} \sin(\pi\alpha n)  \Gamma(1+\alpha n) {x^n\over n!}.$$
Its Mellin transform is an immediate consequence of (1.25), and we have

$$\int_0^\infty  \Phi_\alpha(x) x^{s-1} dx=  { \Gamma (s)\over \Gamma(\alpha s)},\  {\rm Re} (s)= \gamma.\eqno(1.26)$$
Moreover, appealing to  (1.13), Entry 8.4.23.1 in [1, Vol. III] and the Mellin-Parseval identity [2], we will arrive at  the integral representation for $K_{\alpha,\nu}(x)$ in terms of $\Phi_\alpha$ and the modified Bessel function $K_\nu$

$$K_{\alpha,\nu}(x) =  \int_0^\infty  \Phi_\alpha \left(t^\alpha\right) K_\nu\left( {x^{1/\alpha}\over t}\right) {dt\over t}= {1\over \alpha} \int_0^\infty  \Phi_\alpha \left( x t\right) K_\nu\left( t^{- 1/\alpha}\right) {dt\over t}.\eqno(1.27)$$
Further, taking some $\beta \in (0,1)$, we have, accordingly,

$$ \int_0^\infty  \Phi_\beta \left({y\over x^\beta}\right) K_{\alpha,\nu} \left(y^{-\alpha}\right) {dy\over y} = {1\over 2\pi i}  \int_{\gamma-i\infty}^{\gamma+i\infty}   { 2^{\alpha s-2} \Gamma (s)\over \Gamma(\alpha \beta s)} \Gamma\left({\alpha s+\nu\over 2}\right) \Gamma\left({\alpha s-\nu\over 2}\right)  x^{\alpha\beta s} ds$$

$$=  \int_0^\infty  \Phi_{\alpha \beta}  \left(t^\alpha\right) K_{\nu} \left( {x^{-\beta}\over t}\right) {dt\over t}.\eqno(1.28)$$

Let us consider the integral operator of the Mellin convolution type with the kernel $\Phi_\alpha$, namely,

$$(G_\alpha f)(x)\equiv (G_\alpha f(t) )(x)= {1\over x}\int_0^\infty \Phi_\alpha \left({t\over x^\alpha}\right) f(t) dt,\quad x >0,\ 0 < \alpha < 1,\eqno(1.29)$$
where $f$ satisfies  suitable conditions for the existence  of the integral (1.29). Comparing with (1.27), this gives the relation

$$K_{\alpha,\nu}(x) = {x^{-1/\alpha}\over \alpha}   \left(G_\alpha  K_\nu\left( t^{- 1/\alpha}\right) t^{-1} \right)\left(x^{-1/\alpha}\right).\eqno(1.30)$$

{\bf Theorem 3}. {\it Let $ x >0, \ 0< \alpha < 1, \ f \in L_1\left(\mathbb{R}_+; t^{-\gamma} dt\right),\ \gamma \ge 0.$ Then $(G_\alpha f)(x)$ is continuous on the interval $[x_0, \ X_0],\ x_0, X_0 >0$,  having the following Mellin-Barnes representation

$$ (G_\alpha f)(x) = {1\over 2\pi i} \int_{\gamma-i\infty}^{\gamma+i\infty} { \Gamma (s)\over \Gamma(\alpha s)}  f^*(1-s) x^{\alpha s-1} ds,\eqno(1.31)$$
where $f^*(s)$ denotes the Mellin transform $[2], [3]$

$$f^*(s)= \int_0^\infty f(t) t^{s-1} dt.\eqno(1.32)$$
Moreover, the composition of two operators $G_\alpha,\ G_\beta,\ \alpha, \beta \in (0,1)$ obeys the following semigroup property}

$$(G_\beta G_\alpha f)(x) = (G_{\alpha \beta} f)(x).\eqno(1.33)$$

\begin{proof}  Indeed, plugging the right-hand side of (1.25) inside the integral (1.29), we have

$$\left| (G_\alpha f)(x)\right| \le {X_0^{\alpha\gamma} \over 2\pi x_0}\int_0^\infty t^{-\gamma} |f(t)|  dt \int_{\gamma-i\infty}^{\gamma+i\infty} \left|{ \Gamma (s)\over \Gamma(\alpha s)}\right|  |ds|$$

$$=  {X_0^{\alpha\gamma} \over 2\pi x_0}  ||f||_{L_1\left(\mathbb{R}_+; t^{-\gamma} dt\right)} \int_{\gamma-i\infty}^{\gamma+i\infty} \left|{ \Gamma (s)\over \Gamma(\alpha s)}\right|  |ds| < \infty. $$
Consequently, the interchange of the order of integration is allowed and, using (1.32), we end up with (1.31). The continuity of $G_\alpha f$ on the interval $[x_0, \ X_0]$ follows via the uniform convergence of the integral (1.31) on this set.  Further, since the Mellin transform $f^*(1-s)$ is bounded and the gamma-quotient 

$${ \Gamma (s)\over \Gamma(\alpha s)} \in L_1\left(\gamma-i\infty, \gamma+i\infty\right) \cap L_2\left(\gamma-i\infty, \gamma+i\infty\right),\ 0 < \alpha < 1,$$
the composition $G_\beta G_\alpha$ can be calculated directly via (1.31), (1.32) and the Mellin-Parseval equality to  obtain

$$(G_\beta G_\alpha f)(x)  = {1\over x}\int_0^\infty \Phi_\beta \left({t\over x^\beta}\right) (G_\alpha f)(t) dt $$

$$=  {1\over 2\pi i} \int_{\gamma-i\infty}^{\gamma+i\infty} { \Gamma (s)\over \Gamma(\alpha \beta s)}  f^*(1-s) x^{\alpha \beta s-1} ds = (G_{\alpha \beta} f)(x).$$
This proves (1.33) and completes the proof of Theorem 3.

\end{proof}

{\bf Remark 4}. By continuity, when $\alpha \to 1-$, one can define the identity operator $(G_1f) (x)= f(x)$, for instance,  in $L_2\left( \mathbb{R}_+; x^{1-2\gamma} dx\right)$.

{\bf Remark 5}. Considering the convergence of the integral (1.31) in the mean square sense, one can extend the definition of $G_\alpha$ for all $\alpha >0$.  Moreover, the Mellin-Parseval equality yields

$$\int_0^\infty  \left| (G_\alpha f)\left( x\right)\right|^2 x^{1-2\alpha\gamma} dx = {1\over 2\pi \alpha}  \int_{-\infty}^{\infty} \left|{ \Gamma (\gamma+it)\over \Gamma(\alpha (\gamma+it))}  f^*(1-\gamma-it)\right|^2 dt.\eqno(1.34)$$

Consequently, an immediate consequence of (1.31) is the reciprocal formula

$${ \Gamma (s)\over \Gamma(\alpha s)}  f^*(1-s) = \alpha \int_0^\infty  (G_\alpha f)(x) x^{- \alpha s} dx,\  s = \gamma +it,\eqno(1.35)$$
where the integral converges in mean square over $(\gamma-i\infty,\ \gamma+i\infty)$.  Furthermore, the composition $G_\beta G_\alpha,\ \alpha, \beta >0$ can be treated as follows

$$(G_\beta G_\alpha f)(x) = {1\over 2\pi i} \int_{\gamma-i\infty}^{\gamma+i\infty} { \Gamma (s)\over \Gamma(\beta s)}  (G_\alpha f)^*(1-s) x^{\beta s-1} ds$$

$$= {1\over 2\pi \alpha i} \int_{\alpha \gamma-i\infty}^{\alpha \gamma+i\infty} {  \Gamma(s/\alpha) \over \Gamma(\beta s)}  f^*\left(1-{s\over\alpha}\right ) x^{\beta s-1} ds$$

$$= {1\over 2\pi  i} \int_{\gamma-i\infty}^{\gamma+i\infty} {  \Gamma(s) \over \Gamma(\alpha\beta s)}  f^*\left(1-s\right ) x^{\alpha \beta s-1} ds= (G_{\alpha \beta} f)(x).\eqno(1.36)$$
This means that (1.33) is valid for all $\alpha, \beta >0$ under the condition

$${  \Gamma(s) \over \Gamma(\lambda s)}  f^*\left(1-s\right ) \in L_2 \left(\gamma-i\infty,\ \gamma +i\infty\right),\ \lambda = \max(\alpha, \beta, \alpha\beta).\eqno(1.37)$$
 As a consequence of (1.35), the Mellin transform $f^*$ can be represented in the form
 
 $$  f^*(s) =   {\alpha\  \Gamma (\alpha (1-s))\over \Gamma(1-s)} \int_0^\infty  (G_\alpha f)(x) x^{- \alpha (1-s)} dx,\quad  {\rm Re}(s)= 1-\gamma$$
  which gives rise to the relation for the inverse operator $G_\alpha^{-1} = G_{1/\alpha},\ \alpha >0$.

\section{Mapping properties of the transformation (1.11)}

The so-called fractional Kontorovich-Lebedev transform (FKL-transform) (1.11) can be identified with the kernel (1.24), where $\nu=i\tau,\ \tau \in \mathbb{R}.$  We begin with the following result.

{\bf Theorem 4}. {\it Let $\alpha > 0$. If $0< \alpha <1,   f \in L_1\left(\mathbb{R}_+; x^{-\gamma} dx\right),\ \gamma > 0,$  then  the FKL-transform $(1.11)$ is  the Kontorovich-Lebedev transform $(1.2)$ of  the operator $G_{\alpha}$

$$(F_\alpha f)(\tau) = (F  G_\alpha f) (\tau) =  \int_0^\infty K_{i\tau}(x)  (G_\alpha f)(x) dx,\eqno(2.1)$$
where $0< \alpha <1$ and  the integral converges absolutely.  The result is valid as well under   conditions $0 < \alpha \le 1$, 

$$ f^*\left(1-s\right ) \in L_2 \left(\gamma-i\infty,\ \gamma +i\infty\right).\eqno(2.2)$$
Under more strong conditions 

$${\Gamma (s)\over \Gamma(\alpha s)}  f^*(1-s) \in L_1  \left(\gamma-i\infty,\ \gamma +i\infty\right),\eqno(2.3)$$
or

$${\Gamma (s)\over \Gamma(\alpha s)}  f^*(1-s) \in L_2  \left(\gamma-i\infty,\ \gamma +i\infty\right),\eqno(2.4)$$
the result is valid for all $\alpha >0$.}

\begin{proof}  Recalling (1.13), we substitute its right-hand side in (1.11) and interchange the order of integration with the use of (1.32). This implies the formula

$$(F_\alpha f)(\tau)= {1\over 2\pi i}  \int_{\gamma-i\infty}^{\gamma+i\infty}   {2^{\alpha s-2} \Gamma (s)\over \Gamma(\alpha s)} \Gamma\left({\alpha s+ i\tau\over 2}\right) \Gamma\left({\alpha s-i\tau \over 2}\right)   f^*(1-s) ds.\eqno(2.5)$$
The interchange is justified by Fubini's theorem via the estimate

$$\int_0^\infty |f(x)| x^{-\gamma} dx \int_{\gamma-i\infty}^{\gamma+i\infty} \left|  {\Gamma (s)\over \Gamma(\alpha s)} \Gamma\left({\alpha s+i\tau\over 2}\right) \Gamma\left({\alpha s- i\tau\over 2}\right) ds \right|$$

$$\le ||f||_{L_1\left(\mathbb{R}_+; x^{-\gamma} dx\right)}\  B(\alpha\gamma, \alpha\gamma)  \int_{\gamma-i\infty}^{\gamma+i\infty} \left| \Gamma (s) ds\right| < \infty,$$
where $B(a,b)$ is Euler's beta function.  Now, since $0 < \alpha < 1$, the quotient $\Gamma (s)/ \Gamma(\alpha s) \in L_1 \left(\gamma-i\infty,\ \gamma +i\infty\right)$. Moreover,  $f^*\left(1-s\right )$ is bounded. Therefore we employ Entry 8.4.23.1 in [1, Vol. III]  to substitute in (2.5), to interchange the order of integration by virtue of the estimate

$$\int_{\gamma-i\infty}^{\gamma+i\infty}    \left| {\Gamma (s)\over \Gamma(\alpha s)} f^*(1-s) \right| \int_0^\infty \left|K_{i\tau} (x) x^{\alpha s-1} dx ds \right| $$

$$\le ||f||_{L_1\left(\mathbb{R}_+; x^{-\gamma} dx\right)}  \int_{\gamma-i\infty}^{\gamma+i\infty}    \left| {\Gamma (s)\over \Gamma(\alpha s)}  ds \right| \int_0^\infty K_{0} (x) x^{\alpha \gamma-1} dx < \infty\eqno(2.6)$$
and to use (1.31), calculating the inner integral with respect to $x$.   Hence after simple substitutions we end up with representation (2.1).  An alternative case can be realized, recalling condition (2.2), or estimate (2.6) under condition (2.3),  or the Mellin-Parseval equality when $f$ satisfies (2.4). 

\end{proof}
An immediate corollary comes due to Theorem 3 and Remark 5.

{\bf Corollary 4.}  {\it Under conditions of Theorems $3,4$ the following operational relations hold valid

$$(F_{\alpha\beta} f)(\tau) = (F_\alpha\ G_\beta f )(\tau) = (F_\beta\ G_\alpha f )(\tau).\eqno(2.7)$$
Moreover, letting $\beta= 1/\alpha$,  one obtains the decomposition of the Kontorovich-Lebedev transform}

$$(F f)(\tau) = (F_\alpha G_{1/\alpha} f )(\tau) = (F_{1/\alpha} G_\alpha f )(\tau).\eqno(2.8)$$
Finally, in this section, we note that the $L_2$-theory for the Kontorovich-Lebedev transform (see [3, Section 2.3]) permits formula (2.1) when $(G_\alpha f) (x) \in L_2(\mathbb{R}_+;  xdx)$,  and the convergence of the integral is in the mean square sense with respect to the norm in $L_2(\mathbb{R}_+; \\ \tau\sinh(\pi\tau) d\tau)$.  Moreover, the Parseval identity takes place

$$\int_0^\infty \left| (G_\alpha f) (x)\right|^2 x dx = {2\over \pi^2} \int_0^\infty \tau\sinh(\pi\tau) \left| (F_\alpha f)(\tau)\right|^2 d\tau.\eqno(2.9)$$
Due to the parallelogram identity we have for two functions $f,g$

$$\int_0^\infty (G_\alpha f) (x) \overline{(G_\alpha g) (x)}  x dx = {2\over \pi^2} \int_0^\infty \tau\sinh(\pi\tau)  (F_\alpha f)(\tau) \overline{(F_\alpha g) (\tau)} d\tau.\eqno(2.10)$$
More general identities, involving two positive parameters $\alpha, \beta$ read (see (2.7))

$$\int_0^\infty (G_\alpha f) (x) \overline{(G_\beta g) (x)}  x dx = {2\over \pi^2} \int_0^\infty \tau\sinh(\pi\tau)  (F_\alpha f)(\tau) \overline{(F_\beta g) (\tau)} d\tau,\eqno(2.11)$$

$$\int_0^\infty (G_\alpha f) (x) \overline{(G_\alpha G_\beta g) (x)}  x dx = {2\over \pi^2} \int_0^\infty \tau\sinh(\pi\tau)  (F_\alpha f)(\tau) \overline{(F_{\alpha\beta} g) (\tau)} d\tau.\eqno(2.12)$$

{\bf Lemma 1}. {\it Let $0 < \alpha < 1,\ f \in L_1(\mathbb{R}_+).$ Then the operator $G_\alpha : L_1(\mathbb{R}_+) \to L_2 (\mathbb{R}_+ ;  x dx)$ is well-defined by formula $(1.29)$ and bounded, i.e.

$$|| G_\alpha f ||_{L_2 (\mathbb{R}_+ ;  x dx)} \le  \sqrt {{\tan\left(\alpha\pi/2\right)\over 2\pi} } \  ||f||_{L_1(\mathbb{R}_+)}.\eqno(2.13)$$
Besides,  the operator $G_\alpha$ is well-defined by $(1.31)$ in  $L_2(\mathbb{R}_+; x dx)$ and bounded

$$|| G_\alpha f ||_{L_2 (\mathbb{R}_+ ;  x dx)} \le  \sqrt {\alpha} \  ||f||_{L_2(\mathbb{R}_+; xdx)},\eqno(2.14)$$
where the corresponding Mellin transform   $f^*(1-s) \in L_2 (1-i\infty, 1+i\infty)$ and the integral $(1.32)$ converges  in mean with respect to the norm in $L_2 (1-i\infty, 1+i\infty)$.}

\begin{proof} In fact, the Mellin-Barnes representation (1.25) can be extended to $\gamma =0$

$$\Phi_\alpha(x)= {1\over 2\pi } \int_{-\infty}^{\infty} { \Gamma (it)\over \Gamma(\alpha it)} x^{- it} dt,\quad x >0.\eqno(2.15)$$
A a direct consequence, we have by virtue of the Mellin-Parseval identity and the reflection formula for the gamma function 

$$\int_0^\infty \left| \Phi_\alpha(x)\right|^2 {dx\over x} = {\alpha \over \pi^2 } \int_{0}^{\infty} {\sinh(\alpha t) \over \sinh( t)}  dt.\eqno(2.16)$$
The integral on the right-hand side of (2.16) has the value due to Entry 2.4.4.1 in [1, Vol. I]. This gives finally

$$||  \Phi_\alpha ||^2_{L_2 (\mathbb{R}_+ ;  dx/x)} = {\alpha \over 2 \pi }\  \tan\left({\alpha\pi\over 2}\right).\eqno(2.17)$$
Hence, recalling (1.29)  and using  the generalized Minkowski inequality,  we deduce 

$$|| G_\alpha f ||_{L_2 (\mathbb{R}_+ ;  x dx)} = \left( \int_0^\infty \left| \int_0^\infty \Phi_\alpha \left({t\over x^\alpha}\right) f(t) dt \right|^2 {dx\over x}\right)^{1/2} $$

$$\le  \int_0^\infty |f(t)|  \left( \int_0^\infty \left| \Phi_\alpha \left({t\over x^\alpha}\right)\right|^2 {dx\over x}\right)^{1/2} dt = {1\over \sqrt\alpha} ||  \Phi_\alpha ||_{L_2 (\mathbb{R}_+ ;  dx/x)}\  ||f||_{L_1(\mathbb{R}_+)}$$

$$=   \sqrt {{\tan\left(\alpha\pi/2\right)\over 2\pi} } \  ||f||_{L_1(\mathbb{R}_+)}.$$
This proves that the operator $G_\alpha$ is well-defined by (1.29) and its norm in $L_2 (\mathbb{R}_+ ;  x dx)$ satisfies inequality (2.13).  On the other hand, equality (1.34) is valid for $\gamma =0$ under conditions of the lemma, and we have

$$ || G_\alpha f ||^2_{L_2 (\mathbb{R}_+ ;  x dx)} = \int_0^\infty  \left| (G_\alpha f)\left( x\right)\right|^2 x dx = {1\over \pi^2}  \int_{0}^{\infty} { \sinh(\alpha t) \over \sinh( t)} \ \left| f^*\left(1-{it\over \pi}\right)\right|^2 dt$$

$$\le  {\alpha \over 2\pi}  \int_{-\infty}^{\infty}  \left| f^*\left(1+ it \right)\right|^2 dt = \alpha||f||^2_{L_2 (\mathbb{R}_+ ;  x dx)}.$$
Hence we proved (2.14) and completed the proof of Lemma 1. 

\end{proof} 

{\bf Remark 6}. For all $\alpha > 0$ the operator $G_\alpha$ is an isometric  and bounded map from the Banach (Hilbert) space normed by (see (2.4))

$$||f||= {1\over \pi} \left(\int_{0}^{\infty} { \sinh(\alpha t) \over \sinh( t)} \ \left| f^*\left(1-{it\over \pi}\right)\right|^2 dt\right)^{1/2}\eqno(2.18)$$
onto $L_2 (\mathbb{R}_+ ;  x dx)$, having the norm equality (1.34) $|| G_\alpha f ||_{L_2 (\mathbb{R}_+ ;  x dx)} = ||f||.$

\section{Inversion theorem}

The main goal of this section is to invert the fractional Kontorovich-Lebedev transform (1.11).  To do this, we appeal to the $L_2$-theory for the KL-transform [3, Ch. 2] and integral representation (2.1), involving the operator $G_\alpha$.  Then under conditions of Lemma 1, Remark 6  and Theorem 2.4 in [4] the inverse KL-transform gives for almost all $x >0$

$$ (G_\alpha f)(x) = {2\over x\pi^2} {d\over dx} \int_0^\infty \int_0^x \tau\sinh(\pi\tau) K_{i\tau}(y)  (F_\alpha f)(\tau) dy d\tau,\  x, \alpha > 0.\eqno(3.1)$$
Hence, reciprocally, from (1.36) we find the inversion formula in the operator form

$$  f(x) = {2\over \pi^2}  \left(G_{1/\alpha} \bigg[ {1\over t} {d\over dt} \int_0^\infty \int_0^t \tau\sinh(\pi\tau) K_{i\tau}(y)  (F_\alpha f)(\tau) dy d\tau\bigg] \right) (x),\  x, \alpha > 0.\eqno(3.2)$$
On the other hand, rewriting equality (1.31) in the form
$$  x^{-1/\alpha} (G_\alpha f)\left(x^{-1/\alpha}\right) = {1\over 2\pi i} \int_{\gamma-i\infty}^{\gamma+i\infty} { \Gamma (s)\over \Gamma(\alpha s)}  f^*(1-s) x^{-s} ds,\eqno(3.3)$$
we observe via simple  changes of variables that the Mellin transform pair

$$  x^{-1/\alpha} (G_\alpha f)\left(x^{-1/\alpha}\right) \to { \Gamma (s)\over \Gamma(\alpha s)}  f^*(1-s)$$
satisfies the Mellin-Parseval equality (1.34).  But letting $\gamma=0$ (see Remark 6) and   taking equality (2.10),  we work out its left-hand side, assuming a real-valued  function $g$ normed by (2.18)  such that $G_\alpha g \in  L_2 (\mathbb{R}_+ ;   xdx)$.  Hence,  if $G_\alpha f \in  L_2 (\mathbb{R}_+ ;   xdx)$ then the  Mellin-Parseval equality [2, Th. 72] yields (cf. (1.34), (3.3))

$$\int_0^\infty (G_\alpha f) (x) (G_\alpha g) (x)  x dx = {1\over \alpha} \int_0^\infty (G_\alpha f)\left(x^{-1/\alpha}\right) (G_\alpha g)\left(x^{-1/\alpha}\right)  x^{-2/\alpha-1} dx $$

$$= {1\over 2\pi \alpha i}  \int_{-i\infty}^{i\infty}  {\Gamma(s) \Gamma(-s)\over \Gamma(\alpha s) \Gamma(- \alpha s)} f^* (1-s)  g^* (1+s) ds = \int_0^\infty f(y) h_\alpha (y) dy,\eqno(3.4)$$
where assuming $ f  \in  L_2 (\mathbb{R}_+ ;   xdx), h_\alpha \in L_2 (\mathbb{R}_+ ;   dx/x)$, 

$$h_\alpha(y) =  {1\over 2\pi \alpha i}  \int_{-i\infty}^{i\infty}    {\Gamma(s) \Gamma(-s)\over \Gamma(\alpha s) \Gamma(- \alpha s)} g^* (1+s)  y^{-s} ds$$

$$ = {1\over 2\pi  }  \int_{-\infty}^{\infty}   { \sinh(\pi \alpha t) \over \sinh(\pi t)}  g^* (1+it)  y^{-it} dt,\quad y, \alpha  >0.\eqno(3.5)$$
 We observe that the latter integral converges in the mean square sense with respect to the norm in $L_2 (\mathbb{R}_+ ;   dx/x)$.  Moreover, the Mellin-Parseval equality implies

$$||h_\alpha||^2_{L_2 (\mathbb{R}_+ ;   dx/x)} = {1\over 2\pi \alpha^2}  \int_{-\infty}^{\infty}  \left| {\Gamma(it) \Gamma(-it)\over \Gamma(\alpha it) \Gamma(- \alpha it)} g^* (1+it)\right|^2 dt$$

$$= {1\over \pi^2 }  \int_{0}^{\infty}  \left| { \sinh( \alpha t) \over \sinh( t)} g^* \left(1+ {it\over \pi}\right)\right|^2 dt$$

$$\le {\alpha \over \pi^2 }  \int_{0}^{\infty}   { \sinh( \alpha t) \over \sinh( t)} \left|g^* \left(1+ {it\over \pi}\right)\right|^2 dt = \alpha ||g||^2,\quad 0 < \alpha \le 1,\eqno(3.6)$$

$$\infty >  ||h_\alpha||^2_{L_2 (\mathbb{R}_+ ;   dx/x)} \ge  \alpha ||g||^2,\quad \alpha \ge 1.\eqno(3.7)$$
We will seek the function $h_\alpha(y) = y  \mu_{(0,x)} (y)$, where  $\mu_{(0,x)} (y)$ is  a characteristic function of the interval $(0,x),\ x >0$. Then from (3.5) via the reciprocal Mellin transform the function $g_x^*(s)$ can be defined by  the equality
$$  g_x^*(s) = { \alpha \Gamma(\alpha (s-1)) \Gamma(- \alpha (s-1))\over \Gamma(s-1) \Gamma(1-s)} {x^s\over s},\quad {\rm Re}(s)= 1.\eqno(3.8)$$
In the case $\alpha > 1$ one can find the explicit value of the function $g(y)$, inverting the Mellin transform (3.8).  Hence we deduce

$$g_x(y) = {\alpha \over 2\pi i}  \int_{1-i\infty}^{1+i\infty} { \Gamma(\alpha (s-1)) \Gamma(- \alpha (s-1))\over \Gamma(s-1) \Gamma(1-s)} {(x/y)^s\over s}  ds $$

$$=  {1 \over 2\pi}  \int_{-\infty}^{\infty} { \sinh(\pi t)\over \sinh(\pi \alpha t)} {(x/y)^{1+it}\over 1+it}  dt =  {1 \over \pi} \int_0^{x/y} \int_{0}^{\infty} { \sinh(\pi t)\over \sinh(\pi \alpha t)} \ \cos(t \log u ) dt du. \eqno(3.9)$$
The inner integral with respect to $t$ is calculated via Entry 2.5.46.4 in [1, Vol. I], and we derive

$$g_x(y) = {2 \alpha \over \pi}\ \sin(\pi/\alpha)  \int_0^{(x/y)^{1/\alpha}} { u^{\alpha-1} du\over u^{2}+ 2\cos(\pi/\alpha) u + 1}$$

$$=  {2 \alpha \over \pi}\  {\rm Im} \bigg[ \int_0^{(x/y)^{1/\alpha}} { u^{\alpha-1} \over u + e^{-i\pi/\alpha} } du \bigg].\eqno(3.10)$$ 
Finally, appealing to Entry 2.2.6.1 in [1, Vol. I], we write the result in terms of the Gauss hypergeometric function

$$ g_x(y) =  {2 x \over \pi y}   {\rm Im} \bigg[ e^{i\pi/\alpha} {}_2F_1 \left( 1, \alpha; 1+\alpha;\  - \left({x\over y}\right)^{1/\alpha}  e^{i\pi/\alpha} \right).\eqno(3.11)$$
On the other hand, the Mellin-Barnes integral in (3.9) does not change the values if we move the line of integration to the left, i.e

$$ g_x(y) = {1 \over 2\pi i}  \int_{-i\infty}^{i\infty} { \Gamma(\alpha s) \Gamma(1- \alpha s)\over \Gamma(s) \Gamma(1-s)} {(x/y)^{1+s}\over 1+s}  ds $$

$$= {1 \over 2\pi i}  \int_{\gamma -i\infty}^{\gamma+i\infty} { \Gamma(\alpha s) \Gamma(1- \alpha s)\over \Gamma(s) \Gamma(1-s)} {(x/y)^{1+s}\over 1+s}  ds,\eqno(3.12) $$
where $-1/\alpha < \gamma < 0$. Consequently,

$$ g^*_x(1+s) = { \Gamma(\alpha s) \Gamma(1- \alpha s)\over \Gamma(s) \Gamma(1-s)} {x^{1+s}\over 1+s}.\quad {\rm Re}(s)= \gamma.\eqno(3.13)$$
Now, returning to (2.10) and combining with (3.4), we have for $ x >0$
$$ \int_0^x f(y) y dy =   {2\over \pi^2} \int_0^\infty \tau\sinh(\pi\tau)  (F_\alpha f)(\tau) (F_\alpha g_x) (\tau) d\tau.\eqno(3.14)$$
But the transform $(F_\alpha g_x) (\tau)$ can be written by the corresponding formula (2.5) under conditions (2.3), (2.4), and we have after simplification

$$(F_\alpha g_x)(\tau)= {1\over 2\pi i}  \int_{\gamma-i\infty}^{\gamma+i\infty}   {  2^{\alpha s-2} \Gamma(1-\alpha s) \over \Gamma(1-s)}  \Gamma\left({\alpha s+ i\tau\over 2}\right) \Gamma\left({\alpha s-i\tau \over 2}\right)  {x^{1-s}\over 1-s} ds$$

$$= {1\over 2\pi i} \int_0^x  \int_{\gamma-i\infty}^{\gamma+i\infty}   {  2^{\alpha s-2} \Gamma(1-\alpha s) \over \Gamma(1-s)}  \Gamma\left({\alpha s+ i\tau\over 2}\right) \Gamma\left({\alpha s-i\tau \over 2}\right) y^{-s} ds dy ,\eqno(3.15)$$
where $0 < \gamma < 1/\alpha.$ Hence, denoting by 

$$\hat{K}_{\alpha,i\tau} (x) = {1\over 2\pi i}  \int_{\gamma-i\infty}^{\gamma+i\infty}   {  2^{\alpha s-2} \Gamma(1-\alpha s) \over \Gamma(1-s)}  \Gamma\left({\alpha s+ i\tau\over 2}\right) \Gamma\left({\alpha s-i\tau \over 2}\right) x^{-s} ds$$

$$= {1\over 4} \mathop{H_{1,3}^{2,1}}\left({x\over 2^\alpha} ; \  {(0,\alpha) \atop (0,1), (i\tau/2, \alpha/2), (- i\tau/2, \alpha/2)}\right),\eqno(3.16)$$
we establish from (3.14), (3.15)  for almost all $x >0 $ and $\alpha > 1$ the following inversion formula for the fractional Kontorovich-Lebedev transform (1.11)

$$f(x) =   {2\over x \pi^2} {d\over dx} \int_0^\infty \tau\sinh(\pi\tau) \int_0^x \hat{K}_{\alpha,i\tau} (y) dy\   (F_\alpha f)(\tau) d\tau.\eqno(3.17)$$
The results of this section are  summarized by the inversion theorem.

{\bf Theorem 5}. {\it Let $\alpha >0$. The fractional Kontorovich-Lebedev transform $(1.11)$ is a bounded isometric map from the  Hilbert space normed by $(2.18)$ onto $L_2(\mathbb{R}; \tau\sinh(\pi\tau) d\tau)$, where the integral $(1.11)$ converges in the mean square sense with respect to the norm in $L_2(\mathbb{R}; \tau\sinh(\pi\tau) d\tau)$ and the Parseval equality holds valid

$${1\over 2} \int_0^{\infty} { \sinh(\alpha t) \over \sinh( t)} \ \left| f^*\left(1-{it\over \pi}\right)\right|^2 dt = \int_0^\infty \tau\sinh(\pi\tau) \left|(F_\alpha f)(\tau)\right|^2  d\tau,\eqno(3.18)$$
having $f^*(s)$ as the Mellin transform $(1.32)$ whose existence is guaranteed in $L_2(1-i\infty, 1+ i\infty)$ when $\alpha \ge 1$ and by Th. $4$  when $0 < \alpha <1$. Besides, when $\alpha > 1$ the inverse operator for almost all $ x >0$ is given by formula $(3.17)$.}

\begin{proof} Indeed, the boundedness of the operator (1.11) follows immediately from Remark 6 and equality (2.9).  Then the assumption $\alpha > 1$, equalities (2.18),  (3.18) and the Mellin-Parseval identity [2] imply
$$\int_0^\infty |f(x)|^2 x dx = {1\over 2\pi}  \int_{-\infty}^{\infty}   \left| f^*\left(1+ it\right)\right|^2 dt < {1\over 2\pi}  \int_{-\infty}^{\infty}    { \sinh(\pi\alpha t) \over \sinh( \pi t)}\left| f^*\left(1+ it\right)\right|^2 dt = ||f||^2 < \infty.$$
Hence, recalling (3.4), (3.7), we end up with (3.14). This immediately leads to (3.17) after differentiation for almost all $ x >0$.

\end{proof}

{\bf Remark 7}. As we observe from the Mellin-Parseval equality,  the Hilbert space normed by (2.18) is a subspace of  $L_2(\mathbb{R};  x dx)$ when $\alpha \ge 1$,  and the case $\alpha =1$ coincides with the  conventional Kontorovich-Lebedev transform (1.2) and its inversion (cf. (3.1) and Ch. 2 in [4] ) with the Parseval equality

$$\int_0^{\infty}  \left| f (x)\right|^2 x dx =  {2\over \pi^2} \int_0^\infty \tau\sinh(\pi\tau) \left|(F_\alpha f)(\tau)\right|^2  d\tau.$$

Finally, we note that the kernel (3.16) $\hat{K}_{\alpha,i\tau} (x)$ can be expressed in terms of the Poisson integral of the kernel (1.13) $K_{\alpha,i\tau} (x)$.  In fact, similarly as in (3.9) we deduce from (3.16)

$$\hat{K}_{\alpha,i\tau} (x) = {1\over 2\pi i}  \int_0^\infty   K_{\alpha,i\tau} (y)  \int_{\gamma-i\infty}^{\gamma+i\infty}   { \Gamma(1-\alpha s) \Gamma(\alpha s) \over \Gamma(1-s) \Gamma (s)} \left({x\over y}\right)^{-s} {ds dy\over y} $$

$$= {1\over \pi }  \int_0^\infty   K_{\alpha,i\tau} ( xy)  \int_{0}^{\infty}   { \sinh(\pi t) \over \sinh(\pi\alpha t)} \cos(t\log y)  {dt dy\over y} $$

$$= {\sin(\pi/\alpha) \over  \alpha \pi }  \int_0^\infty   {K_{\alpha,i\tau} ( xy) \over y( y^{1/\alpha}+  y^{- 1/\alpha} + 2 \cos(\pi/\alpha) )} dy,$$
i.e. after a simple substitution we establish the integral representation of the kernel $\hat{K}_{\alpha,i\tau} (x)$

$$\hat{K}_{\alpha,i\tau} (x) = {\sin(\pi/\alpha) \over  \pi }  \int_0^\infty   {K_{\alpha,i\tau} \left( xy^\alpha\right) \over  y^{2}+ 2 \cos(\pi/\alpha) y+1}\ dy,\ \alpha > 1, x >0.\eqno(3.19)$$

\bigskip
\centerline{{\bf Funding Declaration}}
\bigskip

\noindent The work was partially supported by CMUP, which is financed by national funds through FCT (Portugal)  under the project with reference UIDB/00144/2025,  https://doi.org/ \\10.54499/UID/00144/2025.

\bigskip
\centerline{{\bf References}}
\bigskip
\baselineskip=12pt
\medskip
\begin{enumerate}

\item[{\bf 1.}\ ] A.P. Prudnikov, Yu.A. Brychkov and O.I. Marichev, {\it Integrals and Series}. Vol. I: {\it Elementary Functions}, Vol. II: {\it Special Functions}, Gordon and Breach, New York and London, 1986, Vol. III : {\it More special functions},  Gordon and Breach, New York and London,  1990.

\item[{\bf 2.}\ ] E.C. Titchmarsh,  {\it An introduction to the theory of Fourier integrals}. New York: Chelsea; 1986.

\item[{\bf 3.}\ ]  S. Yakubovich and Yu. Luchko, The Hypergeometric Approach to Integral Transforms and Convolutions, {\it Kluwer
Academic Publishers, Mathematics and Applications.} Vol.287, 1994. 

\item[{\bf 4.}\ ] S. Yakubovich, {\it Index Transforms}, World Scientific Publishing Company, Singapore, New Jersey, London and
Hong Kong, 1996.

\end{enumerate}

\vspace{5mm}

\noindent S.Yakubovich\\
Department of  Mathematics,\\
Faculty of Sciences,\\
University of Porto,\\
Campo Alegre st., 687\\
4169-007 Porto\\
Portugal\\
E-Mail: syakubov@fc.up.pt\\

\end{document}